\documentclass[11pt]{amsart}
\usepackage{amsmath,amsfonts}
\usepackage{graphicx}

\begin{document}
\title[Edge colorings and Hamiltonian cycles]
{An effective algorithm for the enumeration of edge colorings
and Hamiltonian cycles in cubic graphs}

\date{\today}
\author{V. Ejov}
\address{School of
Mathematics and Statistics, University of South Australia,
Mawson Lakes SA 5095 Australia}
\author{N. Pugacheva}
\address{Institute for Energy Technology, St.Petersburg 197183 Russia}
\author{S. Rossomakhine}
\address{School of
Mathematics and Statistics, University of South Australia,
Mawson Lakes SA 5095 Australia}
\author{P. Zograf}
\address{Steklov Mathematical Institute,
St.Petersburg 191023 Russia}


\begin{abstract}
We propose an effective algorithm  that enumerates (and actually finds)
all 3-edge colorings and Hamiltonian cycles in a cubic graph. The idea is
to make a preliminary run that separates the vertices into two types:
``rigid'' (such that the edges incident to them admit a
unique coloring) and ``soft'' ones (such that the edges incident to them
admit two distinct colorings), and then to perform the coloring.
The computational complexity of this algorithm is on a par with (or even below)
the fastest known algorithms that find a single 3-edge coloring or a
Hamiltonian cycle for a cubic graph.
\end{abstract}

\maketitle

\section{Preliminaries}

Let us recall here some basic facts about the relationship between
3-edge colorings (also called Tait colorings) and Hamiltonian
cycles in cubic graphs; the details can be found in \cite{Tutte}.
By a {\em cubic graph} we understand a connected 3-regular
multi-graph that is allowed to have double edges, but no loops
(obviously, a graph with loops cannot have either edge colorings
or Hamiltonian cycles). Consider a set of three distinct elements
called ``colors'' (say, $\{r,g,b\}$, where $r$ stands for ``red'',
$g$ --- for ``green'', and $b$ --- for ``blue''). A {\em 3-edge
coloring}, or {\em Tait coloring} is an assignment of a color to
every edge such that the edges incident with each vertex have
distinct colors\footnote{Introduced by P.~G.~Tait in 1880 for
planar cubic maps, 3-edge colorings uniquely correspond to
4-colorings of maps.}. Every 3-edge coloring of a cubic graph $G$
gives rise to three distinct 2-factors (that is, 2-regular
spanning subgraphs) of $G$ called {\em Tait cycles}: each Tait
cycle is the union of edges painted in two colors out of the three
(the complement to a Tait cycle is a perfect matching -- the union
of disjoint edges painted in the third color). If a cubic graph
$G$ has a Hamiltonian cycle, then $G$ also admits a 3-edge
coloring that is unique up to a permutation of colors: just paint
the Hamiltonian cycle (which always has even length) in two
intermittent colors, and paint the complement perfect matching in
the remaining third color.

The above connection between 3-edge colorings and Hamiltonian cycles suggests
the following method of enumerating (and actually finding) all the Hamiltonian
cycles in a cubic graph:
\begin{enumerate}
\item Find all 3-edge colorings (Tait colorings) of a given cubic graph up to
permutations of colors;
\item Find all the corresponding 2-factors (Tait cycles);
\item Check for connected 2-factors (Hamiltonian cycles).
\end{enumerate}
This procedure gives the complete list of Hamiltonian cycles in a cubic graph.

For the enumeration of 3-edge colorings we use an exhaustive backtracking
algorithm that works in two runs. During the first run it dynamically separates
all the vertices into two types: ``rigid'' ones (that admit a unique coloring
of edges incident with them), and ``soft'' ones (with exactly two
possibilities of coloring the incident edges); no backtracking is needed on
this stage. The second run is the actual painting of edges: after it successfully
colored the graph or was unable to complete the coloring, it returns to the last
visited soft vertex and tries a different possibility. The details are explained
in the next section.

\section{Description of the algorithm}

First we partition the set of vertices $V(G)$, of the graph $G$
into two disjoint subsets $V(G)= R(G) \cup S(G)$ of rigid ($R$)
and soft ($S$) vertices. We note that this partition is not
canonical. Initially we put $S=R=\emptyset$ and dynamically change
their content. We also introduce a temporary set $U$ of {\em
unidentified} vertices that we already visited, and an ordered
list of colored edges $E$. We label the vertices of $G$ by
integers $\{0,\dots,n-1\},\, n=\#V(G)$. For the vertex with number
$i$ we denote the numbers of adjacent vertices by $n_0^i\leq
n_1^i\leq n_2^i$. An edge connecting $i$ and $j$ we denote by
$[i,j]$. We start at the vertex $0$ and add it to the set $R$ of
rigid vertices. We add the three edges
$[0,n^0_0],\,[0,n^0_1],\,[0,n^0_2]$ incident with it to the list
$E$, and we add their endpoints ${n^0_0, n^0_1, n^0_2}$ to $U$.
Now we check if any of the vertices in $U$ are the endpoints of at
least two edges in $E$. If this is the case, we move all such
vertices from $U$ to $R$, and for every such vertex we also add
the remaining third edge incident to it to the set $E$ (if it is
not already there). We continue the above procedure until there is
no vertex left in $U$ that is an endpoint of at least two edges in
$E$. Now pick the vertex from $U$ with the smallest number, say
$i$, and move it into $S$. Note that $i$ is an endpoint of a
single edge in the current set $E$. Next, we append to $E$ the two
remaining edges incident with $i$, add their endpoints to $U$, and
again check if any of the vertices in $U$ bound at least two edges
in $E$. If they do, such vertices are moved to $R$, and the
missing edges incident with these vertices are added to $E$.
Otherwise, we pick a vertex from $U$ with the smallest number,
move it to $S$ and repeat the procedure until $U$ becomes empty,
or, equivalently, until $E$ coincides with $E(G)$. Since $G$ is
connected, this would mean that $V(G)= R \cup S$. Setting $R(G)=R$
and $S(G)=S$ we obtain the required partition.

The final list of edges $E$ provides the order in which we attempt to
paint the edges of the graph. As above, we start at the vertex 0
and paint the edges
$[0,n^0_0],\,[0,n^0_1],\,[0,n^0_2]$ incident
with it in colors $a,b,c$ respectively. If the next
edge in $E$ has a rigid vertex as its endpoint, then there is a
unique color left for the remaining third edge incident with that
vertex. If the next pair of edges is incident with a soft vertex,
we choose one of the two options
for painting these edges (in case we visit this vertex for the first
time). Continuing this way we may successfully reach the end of
the list $E$ and get a complete edge coloring that we save for our
record. It may as well happen that the procedure ends prematurely
when two edges incident to the same vertex are painted in the same color.
In both cases we return in $E$ to the previous soft vertex that was
visited only once, and start over painting edges in a different way.
At the end we get the list of all
possible Tait colorings and check which of them produce connected
Tait cycles (that is,  Hamiltonian cycles).

Let us now illustrate this algorithm on a simple example.

\begin{figure}
\begin{center}
\includegraphics[scale=0.5]{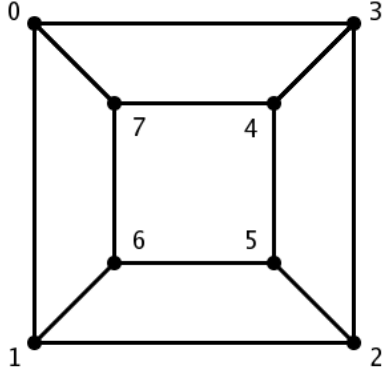}
\caption{1-skeleton of a 3-cube.}
\end{center}
\end{figure}

Consider the graph $G$ on 8 vertices  $\{0,1,2,3,4,5,6,7\}$ shown
on Fig. 1. We start with $S=R=U=E=\emptyset$. The initial vertex
$0$ is a rigid vertex, $R=\{0\}$. Add the edges
$[0,1],[0,3],[0,7]$ to $E$ and the vertices $1,3,7$ to $U$. None of
these vertices is an endpoint of two edges in $E$, so we remove
the vertex 1 from $U$, append it to $S$, and add the edges $[1,2]$
and $[1,6]$ to $E$. Again, none of the vertices in $U$ is an endpoint
of two edges in $E$, so we move the vertex 2 to $S$ and add the edges
$[2,3]$ and $[2,5]$ to $E$. Now the vertex 3 bounds two edges in $E$,
namely, $[0,3]$ and $[2,3]$. We move it to $R$, add the edge $[3,4]$
to $E$ and the vertex 4 to $R$. None of the vertices $\{4,5,6,7\}$ in $U$
bounds at least two edges, so we move the vertex 4 to $S$ and add the
edges $[4,5]$ and $[4,7]$ to $E$. We see that the vertices 5 and 7 become
rigid, and so does the vertex 6. Thus, we get $S(G)=\{1,2,4\},\;
R(G)=\{0,4,5,6,7\}$, and the ordered set of edges is
$E(G)=\{[0,1],[0,3],[0,7],[1,2],[1,6],[2,3],[2,5],[3,4],[4,5],[4,7],[5,6],[6,7]\}$.

Let us start coloring the graph. We paint the edges
$[0,1],[0,3],[0,7]$ in colors $r,g,b$ respectively. There are two
possible ways of coloring the edges incident with the first soft
vertex 1. We first paint $[1,2]$ in color $g$ and $[1,6]$ --- in
color $b$. At the next soft vertex 2 we again have two options,
and choose the first one of them --- paint $[2,3]$ in color $r$,
and $[2,5]$ --- in color $b$. The edge $[3,4]$ incident with the
rigid vertex 3 necessarily have color $b$. At the last soft vertex
4 both possibilities of coloring the edges $[4,5],[4,7]$ lead to
complete edge colorings of $G$. Now we return to the previous soft
point 2 and paint$[2,3]$ in $b$, and $[2,5]$ --- in $r$. This
gives us one more Tait coloring (one of the options cannot be
completed). Finally, returning to the first soft point 1 and
painting the edges $[1,2],[1,6]$ in colors $b,g$ respectively, we
get the last coloring. Thus, there exist 4 distinct Tait colorings
of the graph $G$ (up to permutations of colors) listed in the
following table:

\begin{center}
\bigskip
\begin{tabular}{lllll}
$[0,1]\quad\quad$ &red&red&red&red\\
$[0,3]$&green$\quad$&green$\quad$&green$\quad$&green\\
$[0,7]$&blue&blue&blue&blue\\
$[1,2]$&green&green&green&blue\\
$[1,6]$&blue&blue&blue&green\\
$[2,3]$&red&red&blue&red\\
$[2,5]$&blue&blue&red&green\\
$[3,4]$&blue&blue&red&blue\\
$[4,5]$&red&green&blue&red\\
$[4,7]$&green&red&green&green\\
$[5,6]$&green&red&green&blue\\
$[6,7]$&red&green&red&red
\bigskip
\end{tabular}\end{center}

The first one of these 4 edge colorings has no Hamiltonian cycles
associated with it, whereas the other three ones produce two Hamiltonian
cycles each. They are listed below (we indicate in brackets the
corresponding alternating colors):

\begin{eqnarray*}
  \{0--1--6--5--2--3--4--7--0\}&\quad\mbox{(red--blue)},\\
  \{0--3--4--5--2--1--6--7--0\}&\quad\mbox{(green--blue)},\\
  \{0--1--2--5--6--7--4--3--0\}&\quad\mbox{(red--green)},\\
  \{0--3--2--1--6--5--4--7--0\}&\quad\mbox{(green--blue)},\\
  \{0--1--2--3--4--5--6--7--0\}&\quad\mbox{(red--blue)},\\
  \{0--3--2--5--4--7--6--1--0\}&\quad\mbox{(red--green)}.
\end{eqnarray*}

\section{Computational complexity}

It is clear that the computational complexity of this algorithm is of order
$2^{\#S(G)}.$ Since the set of soft vertices $S(G)$ depends on the ordering
of vertices of $G$, the complexity also depends on this ordering. To give an
upper bound for $\# S(G)$ for a simple graph without double edges
we note that, every time we add two new edges to $E$
incident with a soft vertex, we encounter one of the three possibilities:
\begin{enumerate}
\item Both new endpoints belong to $U$;
\item One new endpoint belongs to $U$ an one to $R$ (the vertex
from $R$ then gives rise to new vertices that are added to $U \cup
R$);
\item Both new endpoints belong to $R$.
\end{enumerate}
In any case, every soft vertex gives rise to at least two new
vertices in $U \cup R$
and the lower bound 2 is attained when the both endpoints of the edges
incident with a soft vertex belong to $U.$ Thus, when the union
$S\cup U\cup R$ becomes equal to $V(G)$ for the first time, we have
the inequality $2\# S \le \# U + \#R.$ In particular, it implies
that $\# S \le n/3$ at this stage.
Let $g$ be the girth of $G$ (that is, the length of the shortest
cycle in $G$). When we reach the stage
$S\cup R\cup U=V(G)$, every new soft vertex gives rise to at
least $g-1$ rigid vertices. This means that no more than $\# U/g$ vertices
will  be added to $S$. Therefore, in the case $(g \ge 4)$, i.e., when
$G$ is triangle free, the number of soft vertices $\#
S\le n/3+ \# U/4 \le n/2,$  as $\# U \le
2n/3$. Thus, the speed of our algoritm is on a par with
the fastest algoritms
that find a single edge coloring or a single Hamiltonian cycle in a cubic graph,
or even better (cf., e.g., \cite{Epp}).
The absence of short cycles makes the algorithm even faster with
complexity bounded from above by $2^{n(g+2)/3g}\approx 2^{n/3}$.
(Note that the presence of double edges does not slow down the algorithm
because at least one of their two common endpoints is rigid.)

It is instructive to compare the above complexity estimate with
the results of \cite{EFLZ}. Let $\lambda_1,\dots,\lambda_n\in [-1,1]$ denote
the eigenvalues of the (normalized) adjacency matrix of a simple cubic
graph $G$. Consider the mean$\mu$ and the variance $\sigma$ of the
exponents $e^{\lambda_i},\;i=1,\dots,n$. For each fixed $n$ the points 
$(\mu,\sigma)$ form clusters called {\em filars} that enjoy a fractal-like
structure. From the results of \cite{EFLZ} combined with the above 
considerations it follows that the closer is the point $(\mu,\sigma)$ to
the origin, the faster works our algorithm for the corresponding graph.

\begin{figure}
\begin{center}
\includegraphics[scale=0.6]{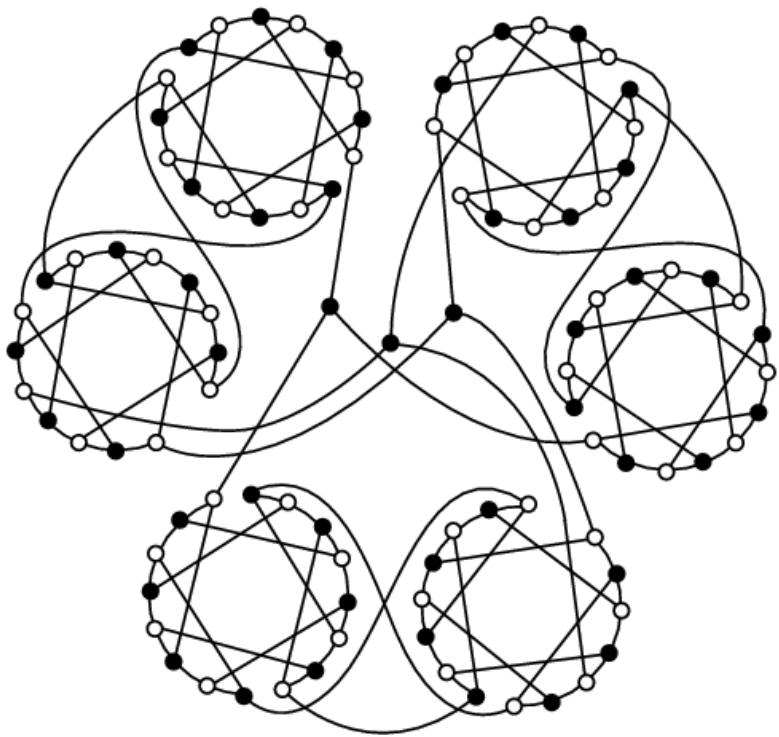}
\caption{Horton graph.}
\end{center}
\end{figure}

This algorithm was implemented in C++ code and compiled on a
Windows x86 machine (Pentium IV 3.40 GHz processor with 1 Gb of
RAM) using CCG GNU Compiler (the program code is given in
Appendix). A good benchmark for testing programs that search for a
Hamiltonian cycle is provided by the Horton graph \cite{Hor},
displayed on Fig. 2.\footnote{This picture is courtesy of:
Weisstein, E. W. ``Horton Graph''. From MathWorld --- a Wolfram
Web Resource; http://mathworld.wolfram.com} It is a cubic
bipartite graph on 96 vertices without Hamiltonian cycles, but
with many ``long'' cycles (that is, cycles of length close to 96).
Some programs choke when they reach such a long cycle, not being
able to transform it into a Hamiltonian one. Our program completed
the search in $5447930319\approx 1.268445\times 10^{32}$ steps (so that
the actual complexity is of order $2^{n/3}$). The process took
6336 sec. of machine time, found 143982592 Tait colorings
and no Hamiltonian cycles.

Our program is an open source program and its ANSI C++ code is available 
at the following address: 
$$\mbox{\small\texttt http://www.unisanet.unisa.edu.au/staff/homepage.asp?Name=Vladimir.Ejov}$$
(we do not present it here because of its length).
The code does not use any platform specific header files,
and with minor modifications it can be compiled with 
essentially any C++ compiler that is not mentally challenged.

{\bf Acknowledgement}
We thank J.~Filar for his interest in this work. The work of VE,
SR was supported, in part, by the Australian Research Council
Discovery grant DP0666632. The work of PZ was partially supported
by the President of Russian Federation grant NSh-U329.2006.1 and
by the Russian Foundation for Basic Research grant 05-01-00899.


\begin{thebibliography}{}
\bibitem{Tutte} Tutte, W. T. {\em Graph Theory}. Cambridge Univ. Press, 1984.
\bibitem{Epp} Eppstein, D. Improved algorithms for 3-coloring,
3-edge-coloring, and constraint satisfaction, 12th ACM-SIAM Symp.
Discrete Algorithms, Washington, 2001, 329--337.
\bibitem{Hor}Bondy, J. A. and Murty, U. S. R., {\em Graph Theory with
Applications}. NY, North Holland, 1976.
\bibitem{EFLZ} Ejov, V. V., Filar, J. A., Lukas, S. K. and Zograf, P.G. 
Clustering of spectra and fractals of regular graphs. 
Preprint math.CO/0610742 (to appear in JMAA).
\end{thebibliography}
\end{document}